\documentclass[12pt]{article}
\usepackage{amstex}
\def\Supp{{\rm Supp\,}}
\def\G{{\rm G}} %the grassmanian
\def\C{{\cal C}}  %the class C
\def\H{{\cal H}}   %the Hausdorff measure 
\newtheorem{lemme}{Lemme}
\newtheorem{theoreme}{Th\'eor\`eme}
\newtheorem{proposition}{Proposition}
\newtheorem{corollaire}{Corollaire}
\newcounter{numeroexemple}
\newenvironment{exemple}
  {\addtocounter{numeroexemple}{1} 
  \begin{trivlist}\item[]\textbf{Exemple \thenumeroexemple}}{\end{trivlist}}
\newenvironment{remarque}
  {\begin{trivlist}\item[]\textbf{Remarque}}{\end{trivlist}}
\newcounter{numerodefinition}
\newenvironment{definition}
  {\addtocounter{numerodefinition}{1}
\begin{trivlist}\item[]\textbf{D\'efinition \thenumerodefinition}}
{\end{trivlist}}
\newenvironment{preuve}{\begin{trivlist}\item[]\textit{Preuve.}}
{\item[] $\square$\end{trivlist}}
\newcommand{\Tan}{{\rm Tan\,}}
\newcommand{\Sing}{{\rm Sing}}
\newcommand{\rang}{{\rm rang}}
\newcommand{\e}{{\rm e}}
\newcommand{\K}{{\cal K}}
\renewcommand{\P}{{\cal P}}
\begin{document}
\title{Conjecture de Globevnik-Stout et 
th\'eor\`eme de Morera pour une cha\^{\i}ne holomorphe}
\author{Tien-Cuong Dinh\footnote{
Institut de Math\'ematiques, UMR 9994 du CNRS,
Universit\'e Pierre et Marie Curie, Tour 46,
Couloir 46--56, 5e \'etage, Case 247,
4, place Jussieu, 75252 Paris Cedex 05, France.
%e--mail dinhtien@math.jussieu.fr
}}
\date{}
\maketitle
{\small
\begin{description} 
\item[{\bf R\'esum\'e.}] Soient $D\subset\subset\mathbb{C}^n$ 
une vari\'et\'e complexe
  de dimension $p\geq 2$ \`a bord
  $\C^2$ dans $\mathbb{C}^n$, 
  $f$ une fonction $\C^1$ sur $bD$ et $V$ une famille
  suffisamment grande et g\'en\'erique de $(n-p+1)$-plans complexes. 
  Supposons que pour $\nu\in V$,
  aucune composante connexe de $bD\cap \mathbb{C}^{n-p+1}_\nu$ n'est
  presque r\'eelle analytique   et que $f$ se prolonge holomorphiquement dans
  $D\cap\mathbb{C}^{n-p+1}_\nu$.  Alors $f$ se prolonge
  en une fonction holomorphe dans $D$. Dans un cas particulier, 
  ce r\'esultat donne une
  r\'eponse partielle  \`a une conjecture de Globevnik-Stout.
En g\'en\'eralisant le th\'eor\`eme de Harvey-Lawson, nous
  d\'emontrons un th\'eor\`eme
  du type de Morera pour le probl\`eme du bord dans $\mathbb{C}^n$ qui
  r\'epond \`a un probl\`eme de
  Dolbeault-Henkin.
\item[{\bf Abstract.}] Let  $D\subset\subset\mathbb{C}^n$ 
be a complex manifold of  dimension
$p\geq 2$ with  $\C^2$ boundary in  $\mathbb{C}^n$.
Let  $f$ be a $\C^1$ function on $bD$ and  $V$ a generic and large enough
family of complex $(n-p+1)$-planes.
Let suppose that for  $\nu\in V$, no connected component of
$bD\cap \mathbb{C}^{n-p+1}_\nu$ is "almost" real analytic and that $f$
extends holomorphically in $D\cap\mathbb{C}^{n-p+1}_\nu$.
Then  $f$ extend as a holomorphic function in $D$. In a special case,
this result gives a
partial answer to a conjecture of Globevnik-Stout.
By generalizing the theorem of Harvey-Lawson, we prove a Morera type theorem
for the boundary problem in $\mathbb{C}^n$ which
answer to a problem asked by
Dolbeault and Henkin.
\end{description}
}
\section{Introduction} Soit $\Gamma$ une hypersurface r\'eelle d'une
  vari\'et\'e complexe $X$ de dimension $n$. 
Une fonction $f$ localement
int\'egrable sur $\Gamma$ est appel\'ee {\it fonction CR} 
si pour toute $(n,n-2)$-forme
 $\psi$ de classe $\C^\infty$ \`a support compact 
 dans $U$, on a $\int_{\Gamma}
 f\overline \partial\psi=0$, o\`u $\Gamma$ poss\`ede l'orientation
 induite par celle de $X$. Si $f\in\C^1$, cette condition est
 \'equivalente \`a condition 
\begin{eqnarray}
\overline Lf(z)=0 & & \mbox{ pour tout vecteur tangent holomorphe} \nonumber\\
& & \mbox{ complexe } 
L\in\mathbb{C}\otimes_{\mathbb{R}}
\Tan_\mathbb{C}(\Gamma,z) 
\end{eqnarray}
pour tout $z\in \Gamma$, o\`u $\Tan_\mathbb{C}(\Gamma,z)$ est 
le sous-espace  tangent complexe maximale de l'espace tangent de $\Gamma$ en
$z$. \\
D'apr\`es le th\'eor\`eme de Bochner, toute
  fonction continue, CR sur le bord lisse d'un domaine simplement
  connexe $D$ de $\mathbb{C}^n$ avec $n\geq 2$ se prolonge contin\^ument
  en une fonction holomorphe de $D$. En g\'en\'eralisation ce
  th\'eor\`eme, Harvey et Lawson ont d\'emontr\'e que toute
  vari\'et\'e compacte $\Gamma$ de classe $\C^1$ \`a singularit\'e
  n\'egligeable, de dimension $2p-1$, orient\'ee et {\it maximalement
    complexe} (c.-\`a-d. le plan tangent de $\Gamma$ en tout point
  contient un $\mathbb{C}^{p-1}$) dans
  $\mathbb{C}^n$ borde une vari\'et\'e complexe $D$ (ou une
  $p$-cha\^{\i}ne holomorphe si $\Gamma$ n'est pas irr\'eductible)
  pour tout $p\geq 2$. Si
  $\Gamma$ est de classe $\C^k$ avec $k\geq 1$, $D$ est $\C^k$ en tout
  point du
  bord sauf sur un compact de volume $(2p-1)$-dimensionnel nul
  \cite{Harvey, HarveyLawson}. On en d\'eduit que
  toute fonction $f$ CR, de classe $\C^1$ sur $bD$ se prolonge en
  une fonction holomorphe, born\'ee dans $D\setminus \Sing D$. En
  effet, le graphe de $f$ est maximalement complexe, elle borde le
  graphe d'une fonction holomorphe sur $D\setminus\Sing D$.\\
Dans cet article, nous allons \'etudier le probl\`eme d'extension
  holomorphe de
  fonctions  sans hypoth\`ese ``$f$ est CR'' et
  \'egalement le probl\`eme du bord sans hypoth\`ese ``$\Gamma$ est
  maximalement complexe''.
\bigskip\\
Soient $D$ un domaine born\'e \`a bord $\C^2$ dans $\mathbb{C}^n$ et
$f$ une fonction continue sur $bD$. D'apr\`es  
Agranovski-Semenov, Valski, Rudin, Stout,
si $f$ se prolonge holomorphiquement dans $D\cap\mathbb{C}_\nu$ pour
$\H^{4(n-1)}$-presque toute droite complexe $\mathbb{C}_\nu$, alors
elle se prolonge holomorphiquement dans $D$
\cite{AgranovskiSemenov,Rudin,Stout1}. Plus g\'en\'eralement, Globevnik
et Stout ont d\'emontr\'e que si sur
$bD\cap\mathbb{C}_\nu$, la fonction $f$
v\'erifie la condition de moment faible, qui s'appelle {\it condition
  de Morera} (c.-\`a-d. $\int_{bD\cap\mathbb{C}_nu}fdz_1=0$) pour
$\H^{4(n-1)}$-presque tout $\nu$, alors elle est CR et elle se
prolonge holomorphiquement dans $D$ si $bD$ est connexe 
\cite{GlobevnikStout}. Ce
th\'eor\`eme reste valable si on consid\`ere uniquement les droites
dont la direction appartient \`a un ouvert non vide de $\G(1,n)$. 
Dans la d\'emonstration, Globevnik et Stout ont consid\'er\'e 
les fonctions dont les lignes de niveaux sont les hyperplans complexes
parall\`eles. L'espace
engendr\'e par telles fonctions est dense dans $\C^\infty$. 
Cette propri\'et\'e r\'eduit le
probl\`eme \`a une v\'erification de l'\'egalit\'e 
$\int_\Gamma f\overline\partial
\psi=0$ pour $\psi=A dz_1\wedge\ldots\wedge dz_n\wedge d\overline
z_1\wedge\ldots\wedge\widehat{d\overline
  z_i}\wedge\ldots\wedge\widehat{d\overline z_j}\wedge\ldots\wedge
d\overline z_n$, o\`u $A$ est une fonction dont les lignes de niveaux
sont des hyperplans complexes parall\`eles.
Finalement, cette \'egalit\'e est un corollaire du th\'eor\`eme de
Fubini et de la condition de Morera.
Globevnik et Stout ont conjectur\'e que si $D$ est convexe,
$\Omega\subset\subset D$ est un sous-domaine convexe et si $f$ se
prolonge holomorphiquement dans $D\cap\mathbb{C}_\nu$ pour toute
droite tangente \`a $b\Omega$, alors $f$ se prolonge holomorphiquement
dans $D$. Notre premier r\'esultat donne une r\'eponse partielle
\`a cette conjecture.\\ 
Soient $Y$ un compact {\it $(n-p+1)$-lin\'eairement convexe} de
  $\mathbb{CP}^n$ (c.-\`a.-d. $\mathbb{CP}^n\setminus Y$ est la
  r\'eunion d'une famille continue de $(n-p+1)$-plans projectifs), 
  $D$ une vari\'et\'e \`a bord $\C^2$ dans $\mathbb{C}^n\setminus Y$,
  born\'ee dans $\mathbb{C}^n$, $f$
  une fonction $\C^1$ d\'efinie sur $bD$ et $V\subset\G(n-p+2,n+1)$
  une famille de $(n-p+1)$-plans complexes de $\mathbb{C}^n\setminus
  Y$ qui coupent $bD$
  transversalement. Supposons que pour tout $\nu\in V$,  
  $f$ se prolonge holomorphiquement
  dans $D\cap\mathbb{C}^{n-p+1}_\nu$ et 
  qu'aucune composante connexe de $bD\cap\mathbb{C}^{n-p+1}_\nu$ n'est
  {\it presque r\'eelle analytique} 
  (c.-\`a-d. $\mathbb{C}_\nu\cap bD$ n'est pas 
  r\'eelle analytique dehors d'un compact de longueur
  $\H^1$ nulle). Nous allons
  d\'emontrer que si $V$ est suffisament grande et g\'en\'erique (par
  exemple une vari\'et\'e de codimension r\'eelle $1$ telle que
  $\bigcup_V\mathbb{C}^{n-p+1}_\nu$ recouvre $bD$) alors $f$ se prolonge
  holomorphiquement dans $D\setminus\Sing D$. Ce th\'eor\`eme n'est
  plus valable si l'on supprime la condition 
  ``aucune composante connexe de $bD\cap\mathbb{C}^{n-p+1}_\nu$ n'est
   presque r\'eelle analytique''. Nos contre-exemples sont donn\'es
  au cas o\`u $bD$ est Levi plat et $Y$ n'est pas un compact de
  $\mathbb{C}^n$. L'id\'ee de la d\'emonstration (par exemple pour
  $n=p=2$, $D$ convexe et $V$ une hypersurface r\'eelle de $\G(2,3)$)
  est la suivante. On consid\`ere $\nu=(\zeta,\eta)\in V$ un point
  g\'en\'erique,  
  $\mathbb{C}_\nu= \{z_2=\zeta+\eta z_1\}$ et $H$ le plan tangent
  de $V$ en $\nu$. Supposons, par exemple, que
  $H\otimes_{\mathbb{R}}\mathbb{C}$ est engendr\'e par
    $\partial/\partial \zeta +\partial/\partial \overline\zeta$,
    $\partial/\partial \eta$ et  $\partial/\partial\overline
    \eta$. Par hypoth\`ese, 
    les d\'eriv\'ees de $R(fPdz_1)=\int_{bD\cap\mathbb{C}_\nu}fPdz_1$ par
    les vecteurs pr\'ec\'edents sont nulles en $\nu$ pour tout
    polyn\^ome $P$ en $z$. D'autre part, 
\begin{eqnarray*}
\frac{\partial R(fPdz_1)}{\partial \eta} & = & 
\frac{\partial R(fz_1Pdz_1)}{\partial \zeta}\\
\frac{\partial R(fPdz_1)}{\partial \overline\zeta} & = & 
\int_{bD\cap\mathbb{C}_\nu}\frac{\partial f}{\partial\overline z_2}
\frac{\partial \overline z_2}{\partial\overline \zeta} P dz_1\\ 
\frac{\partial R(fPdz_1)}{\partial \overline\eta} & = & 
\int_{bD\cap\mathbb{C}_\nu}\frac{\partial f}{\partial\overline z_2}
\frac{\partial \overline z_2}{\partial\overline \zeta}\overline z_1 P dz_1\\ 
\end{eqnarray*}
Ces \'egalit\'es sont induites par le lemme de Darboux \cite{Henkin}
dans le cas o\`u $bD$ est, au voisinage de
$bD\cap\mathbb{C}_\nu$,
feuill\'et\'e par les droites complexes. 
Pour le cas g\'en\'eral, nous
approximerons $bD$ par ses droites complexes tangentes en
$bD\cap\mathbb{C}_\nu$. Nous d\'eduisons des \'egalit\'es pr\'ec\'edentes,
appliqu\'ees pour $P$ ou pour $P:=z_1P$, que $\frac{\partial
f}{\partial\overline z_1}\frac{\partial \overline
z_2}{\partial\overline \zeta}\overline z_1$ et  $\frac{\partial
f}{\partial\overline z_1}\frac{\partial \overline
z_2}{\partial\overline \zeta}z_1$ v\'erifient la condition
des moments. Par cons\'equent, elles se prolongent holomorphiquement
dans $D\cap\mathbb{C}_\nu$. 
Si $bD\cap \mathbb{C}_\nu$ n'est pas presque r\'eelle
analytique, la fonction $\overline z_1$ ne se prolonge pas
m\'eromorphiquement dans $D\cap\mathbb{C}_\nu$. Les deux fonctions
pr\'ec\'edentes sont donc nulles. D'o\`u $f$ v\'erifie (1) sur
$bD\cap\mathbb{C}_\nu$. 
\bigskip\\
Un ferm\'e $K$ d'une vari\'et\'e $X$ est
  appel\'e {\it $k$-g\'eom\'etriquement rectifiable} s'il est
  localement $(\H^k,k)$-rectifiable et s'il admet un $k$-plan tangent
  r\'eel g\'eom\'etrique 
  en tout point sauf sur un ensemble de {\it mesure de
  Hausdorff $k$-dimensionnelle} $\H^k$ nulle. Cette notion permet de
  g\'en\'eraliser le th\'eor\`eme de Harvey-Lawson pour tout courant
  rectifiable, ferm\'e, maximalement complexe et \`a support compact
  g\'eom\'etriquement $(2p-1)$-rectifiable dans $\mathbb{C}^n$ 
  \cite{Dinh2}.\\ 
Soit $\Gamma$ un
  courant rectifiable, ferm\'e, de dimension $2p-1$ et \`a support
  compact g\'eom\'etriquement $(2p-1)$-rectifiable tel que le courant
  d'intersection $\Gamma\cap\mathbb{C}^{n-p+1}_\nu$ v\'erifie {\it la
  condition de Morera} ($(\Gamma\cap\mathbb{C}^{n-p+1}_\nu,z_idz_j)=0$
  pour tous $i,j$)
  pour tout $\nu\in \G(n-p+2,n+1)$
  g\'en\'erique. Alors $\Gamma$ borde une $p$-cha\^{\i}ne holomorphe
  au sens des courants. Ce th\'eor\`eme donne la r\'eponse positive
  \`a un probl\`eme de Dolbeault-Henkin \cite[probl\`eme
  II]{DolbeaultHenkin} et g\'en\'eralise les r\'esultats de
  d'Agranovski-Semenov, Rudin, Valski, Stout, Globevnik
  \cite{Stout1,GlobevnikStout}. Il g\'en\'eralise \'egalement
  le th\'eor\`eme de Harvey-Lawson. En effet,
gr\^ace \`a la formule de Cauchy, on peut montrer que si $\Gamma$ est
maximalement complexe, 
$\Gamma\cap\mathbb{C}^{n-p+1}_\nu$ v\'erifie la condition des
moments pour un sous-espace $\mathbb{C}^{n-p+1}_\nu$ g\'en\'erique
de $\mathbb{C}^n$. La solution locale du probl\`eme du bord 
au voisinage d'un point de
  convexit\'e (Lewy) 
permet de recoller les solutions du probl\`eme du bord sur les tranches
$\mathbb{C}^{n-p+1}_\nu$ en une $p$-cha\^{\i}ne holomorphe. Dans la
d\'emonstration du th\'eor\`eme de Morera, 
nous utilisons l'id\'ee de Globevnik-Stout sur la densit\'e
d'un certain sous-espace de fonctions $\C^\infty$ de $\mathbb{C}^n$.\\ 
Nous d\'emontrons \'egalement un th\'eor\`eme analogue \`a la
conjecture de\break
Globevnik-Stout pour le probl\`eme du bord dans $\mathbb{C}^n$.
\section{Probl\`eme d'extension holomorphe}
Soient $D$ un domaine de Jordan \`a bord rectifiable de
$\mathbb{C}$ et $f$ une fonction continue sur $bD$. Alors $f$ se prolonge
contin\^ument 
en une fonction holomorphe dans $D$ si et seulement si $f$ v\'erifie
{\it la condition des moments}, c.-\`a-d. $\int_{bD} f\varphi=0$ pour toute
$(1,0)$-forme polynomiale $\varphi$ sur $\mathbb{C}$. Plus
g\'en\'eralement, soient $D$ une surface de Riemann polynomialement
convexe de $\mathbb{C}^n$
\`a bord  g\'eom\'etriquement 1-rectifiable et $f\in{\cal
  C}^0(bD)$ v\'erifiant {\it la condition des moments} dans $\mathbb{C}^n$,
c.-\`a-d. $(d[D],f\varphi)=0$ pour toute $(1,0)$-forme polynomiale
$\varphi$ dans $\mathbb{C}^n$, o\`u $[D]$ est le courant de
bidimension $(1,1)$ d\'efini par l'int\'egration sur $D$. Alors $f$
se prolonge en une fonction holomorphe dans $D$ au sens faible des
courants. Ceci est un corollaire
d'un r\'esultat de \cite{Dinh3} sur l'extension des mesures
orthogonales en $(1,0)$-forme holomorphe 
(appliqu\'e aux mesures $f\varphi$ support\'ees par
$bD$). Ce r\'esultat de \cite{Dinh3} 
g\'en\'eralise des r\'esultats de Wermer et de Henkin, qui sont valables au
cas d'une courbe r\'eelle analytique \cite{Wermer} ou 
$\C^2$ par morceaux \cite{Henkin}. 
\begin{definition} (\cite{Henkin})
Une fonction
(ou une $(1,0)$-forme) m\'eromorphe $f$ sur $D$ est appel\'ee {\it
  holomorphe} 
si le courant d'int\'egration $f\wedge [D]$ est
$\overline\partial$ ferm\'e dans $\mathbb{C}^n\setminus bD$. Si $D$
est lisse, cette d\'efinition donne les fonctions et les formes
holomorphes habituelles.
\end{definition}
Soit $\gamma\subset\mathbb{C}^n$ une courbe r\'eelle, ferm\'ee, de
classe $\C^2$ bordant une surface de Riemann $S'$ (\'eventuellement
singuli\`ere et r\'eductible). Alors l'enveloppe polynomiale
$\widehat\gamma$  de
$\gamma$ est la r\'eunion de $\gamma$ avec une surface de Riemann $S$
contenant la surface pr\'ec\'edente \cite{Stolzenberg}. 
On note $S_j$ avec $j\in J$ les
surfaces de Riemann irr\'eductibles de
$S$. On appelle $\gamma_j$ la courbe
ferm\'ee minimale de $\gamma$ telle que $\widehat\gamma_j\supset 
S_j$. D'apr\`es le th\'eor\`eme d'unicit\'e, 
la surface $S$ est lisse jusqu'au bord 
sauf sur un ensemble
d\'enombrable de singularit\'e de $S$ et sur un compact de mesure $\H^1$
nulle (\'eventuellement vide) 
de $\gamma$. Toute mesure orthogonale $\mu=hdz_m$ se
prolonge en une $(1,0)$-forme holomorphe
$\varphi dz_m$ dans $S$ au
sens faible des courants. De plus, si $\gamma_j\not\subset bS_{j'}$ pour tout
$j'\not =j$, alors pour $\H^1$-presque tout point
$x\in\gamma_j$ la fonction $\varphi(z)$ tend vers $h(x)$ quand $z$ tend
vers $x$ le long des arcs non tangentiels \`a $\gamma_j$ \cite{Dinh3}.
Si $h$ est une fonction continue, il existe un compact
$K\subset\gamma_j$ de longueur z\'ero ($\H^1(K)=0$) tel que le
prolongement de $h$ soit continu en tout point de $\gamma_j\setminus
K$. 
On note $\K_j$ 
l'anneau des fonctions $h$ continues sur $\gamma_j$ 
telles que $h|_{\gamma_j\setminus K_h}$ se prolonge en une fonction
m\'eromorphe sur $S_j$, 
continue en tout point de $\gamma_j\setminus K_h$ pour un
certain compact $K_h\subset\gamma_j$ de longueur 0. En particulier,
dans un voisinage de $\gamma_j\setminus K_h$ l'extension de $h$ n'a
pas de p\^ole.
La courbe $\gamma$ est appel\'ee {\it g\'en\'erique} si pour tout
$j\in J$ 
les fonctions $\overline z_m$ sont $\K_j$-alg\'ebriquements
ind\'ependantes sur $\gamma_j$,
c.-\`a-d. pour tout polyn\^ome $P$ de $\K_j[x_1,x_2,\ldots,
x_n]$, on a  $P(\overline z)|_{\gamma_j}\not\equiv 0$. En particulier,
si $\gamma$ est g\'en\'erique, elle n'est incluse dans aucune
hypersurface alg\'ebrique de $\mathbb{C}^n$ et l'application naturelle
de l'anneau des polyn\^omes $\P$ dans $\K_j$ est injective. \\
On note $\G(q+1,n+1)$ la grassmannienne qui param\`etre l'ensemble des
$q$-plans projectifs de $\mathbb{CP}^n$. On a
$\dim\G(q+1,n+1)=(q+1)(n-q)$. L'ensemble
$\G^*(q+1,n+1)$ des $q$-plans affines de $\mathbb{C}^n$ sera
consid\'er\'e comme un ouvert de Zariski de $\G(q+1,n+1)$. Soit
$Y\subset\mathbb{CP}^n$ un compact  $q$-lin\'eairement convexe. On pose
$\G^*_Y(q+1,n+1)=\{\nu\in\G^*(q+1,n+1):\ \mathbb{C}^q_\nu\cap
Y=\emptyset\}$. Soit $D$
une vari\'et\'e complexe \`a bord $\C^2$ de dimension $p\geq 2$ 
dans $\mathbb{C}^n\setminus Y$, born\'ee dans $\mathbb{C}^n$,
c.-\`a-d. $D$ est un sous-ensemble analytique de dimension pure $p$ de
$\mathbb{C}^n\setminus bD\cup Y$. Alors $bD$ admet une orientation induite par
celle de $D$ ainsi que $bD\cap \mathbb{C}^{q}_\nu$ si cette
intersection est transversale. On note
$\G_D(q+1,n+1)$ le sous-ensemble de $\G^*_Y(q+1,n+1)$ 
des $q$-plans  qui coupent $bD$
transversalement. D'apr\`es le th\'eor\`eme de Sard, cet ensemble
est un ouvert dense de $\G^*_Y(q+1,n+1)$ et son compl\'ementaire
est de volume $0$. 
\begin{theoreme} Soient $Y$ un compact $(n-p+1)$-lin\'eairement 
   convexe de $\mathbb{CP}^n$,
   $D$ une vari\'et\'e de dimension $p\geq 2$  
  de $\mathbb{C}^n\setminus Y$, born\'e dans
  $\mathbb{C}^n$, \`a bord $\C^2$ et $f$ une fonction $\C^1$ d\'efinie
  sur $bD$. Soit $V\subset \G_D(n-p+2,n+1)$ tel que $bD\cap\bigcup_{\nu\in
  V_1}\mathbb{C}_\nu$ soit dense dans $bD$, 
o\`u $V_1$ est l'ensemble de $\nu\in V$
  v\'erifiant 
\begin{enumerate}
\item Le c\^one tangent de $V$ en $\nu$ contient $2(p-1)(n-p+2)-1$ 
vecteurs  r\'eels ind\'ependants. 
\item Aucune composante connexe de $\mathbb{C}^{n-p+1}_\nu\cap 
bD$ n'est presque
  r\'eelle analytique.
\end{enumerate}
Supposons que $f$ se prolonge contin\^ument en une fonction holomorphe
dans $D\cap\mathbb{C}^{n-p+1}_\nu$ pour tout $\nu\in V$. Alors $f$ se prolonge
contin\^ument en une fonction holomorphe dans $D\setminus\Sing D$,
localement born\'ee dans $\overline D\setminus Y$.
\end{theoreme}
\begin{remarque} Dans le cas g\'en\'eral ($Y$ n'est pas un compact de
  $\mathbb{C}^n$) la condition 2 n'est pas supprimable ({\it voir} les
  exemples 1 et 2).
\end{remarque}
Soient $\nu_0$ un point de $\G^*(n-p+2,n+1)$, $U$ un
  voisinage de $\mathbb{C}^{n-p+1}_{\nu_0}$, $B\subset U$ une
  sous-vari\'et\'e r\'eelle, orient\'ee, 
  maximalement complexe de dimension $2p-1$, 
  born\'ee dans $\mathbb{C}^n$, 
  \`a bord $\C^2$, $f$ une fonction $\C^1$ sur $B$ et $\varphi$ une
  $1$-forme d\'efinie sur $B$. Supposons que $\mathbb{C}^{n-p+1}_\nu$
  coupe $B$ transversalement. On appelle (en g\'en\'eralisation 
  la transformation d'Abel-Radon d\'efinie dans \cite{Henkin}) 
{\it la transformation
  d'Abel-Radon} de $[B]\wedge f\varphi$ l'int\'egrale $\int_{B\cap
  \mathbb{C}^{n-p+1}_\nu} f\varphi$. Cette int\'egrale d\'efinit une fonction
  au voisinage de $\nu_0$, continue si $f$ est continue et
  de classe $\C^1$ si $f$ l'est.
\begin{proposition} 
  Supposons que $\gamma=\mathbb{C}^{n-p+1}_{\nu_0}\cap B$
  borde une
  surface de Riemann $S'$. Supposons que pour toute $(1,0)$-forme
  polynomiale $\varphi$ de $\mathbb{C}^n$, la d\'eriv\'ee de la
  transformation d'Abel-Radon 
  $R(f\varphi)$ par rapport \`a $\nu$ s'annulle en $\nu_0$ sur un $u$-plan
  r\'eel g\'en\'erique $H\subset\Tan(\G^*(n-p+2,n+1),\nu_0)$ 
  ind\'ependant de $\varphi$ avec $u\geq (p-1)(n-p+2)+1$. L'une des conditions
  suivantes sera satisfaite
\begin{enumerate} 
\item $\overline L
  f(z)=0$ pour tout $z\in \gamma$ et
  pour tout vecteur tangent complexe holomorphe  
  $L\in\mathbb{C}\otimes_{\mathbb{R}} \Tan_\mathbb{C}(B,z)$.
\item Il existe une courbe $\gamma_j$ avec $j\in  J$
   telle que  
  les conjug\'ees des coordonn\'ees de $\mathbb{C}^{n-p+1}_{\nu_0}$ 
  sont  $\K_j$-alg\'ebriquement
  d\'ependantes.
\end{enumerate}
\end{proposition}
\begin{remarque} Les courbes $\gamma_j$ et les anneaux $\K_j$ sont
   d\'efinis au dessus.
La notion de $u$-plan g\'en\'erique sera d\'efinie plus tard. Si
$u=2(p-1)(n-p+2)-1$ ou $2(p-1)(n-p+2)$, 
tout $u$-plan est g\'en\'erique dans ce sens ({\it voir} le lemme
2).
\end{remarque}
\begin{preuve} Dans un ouvert $\mathbb{C}^{(p-1)(n-p+2)}$ de Zariski de
    $\G(n-p+2,n+1)$, 
on peut identifier $\nu$ \`a une matrice $(\zeta,\eta)$ de taille 
$(p-1)\times (n-p+2)$, o\`u 
$$\zeta=\left(
\begin{array}{c}
\zeta_{n-p+2}\\
\zeta_{n-p+3}\\
\vdots\\
\zeta_n
\end{array}
\right)
\mbox{ et } 
\eta=\left(
\begin{array}{cccc}
\eta^1_{n-p+2} & \eta^2_{n-p+2} & \cdots & \eta^{n-p+1}_{n-p+2}\\
\eta^1_{n-p+3} & \eta^2_{n-p+3} & \cdots & \eta^{n-p+1}_{n-p+3}\\
\vdots & \vdots & \ddots & \vdots\\
\eta^1_n & \eta^2_{n} & \cdots & \eta^{n-p+1}_{n}
\end{array}
\right)$$
et le $(n-p+1)$-plan complexe $\mathbb{C}^{n-p+1}_{\nu}$ est d\'efini
par les $(p-1)$ \'equations suivantes
\begin{eqnarray}
z_m & = & \zeta_m+\eta^1_m z_1+\eta^2_m z_2+\cdots + \eta^{n-p+1}_m
z_{n-p+1}\nonumber\\
& &  \mbox{ pour } m=n-p+2,\ldots,n
\end{eqnarray}
Consid\'erons le cas o\`u $B$ est une r\'eunion de $(p-1)$-plans
complexes, c.-\`a-d. $B=\bigcup_{a\in\Upsilon}\mathbb{C}^{p-1}_a$, o\`u
$\Upsilon$ est une courbe r\'eelle, orient\'ee, ferm\'ee (\'eventuellement
r\'eductible) de $\G^*(p,n+1)$. Sans perdre en g\'en\'eralit\'e, on
suppose que pour tout $a\in\Upsilon$, $z'':=(z_{n-p+2},z_{n-p+3},\ldots,z_n)$
est un syst\`eme des coordonn\'ees de $\mathbb{C}^{p-1}_a$ et
$z':=(z_1,z_2,\ldots,z_{n-p+1})$ est un syst\`eme des coordonn\'ees 
de $\mathbb{C}^{n-p+1}_{\nu_0}$. Dans ce
cas, on peut \'ecrire la fonction $f$ sous la forme d'une fonction
$g(a,z'')$ \`a
variables $a$ et $z''$. Pour tout $a\in\Upsilon$ et tout $\nu$ pr\`es de
$\nu_0$, on note $Z:=(Z_1,Z_2,\ldots,Z_n)$ les
coordonn\'ees de point d'intersection de $ \mathbb{C}^{p-1}_a$ et de
$\mathbb{C}^{n-p+1}_{\nu}$, qui d\'ependent de $a$ et de $\nu$.\\
Les \'equations (2) impliquent 
\begin{eqnarray}
\frac{\partial Z_m}{\partial \eta^k_q} & = & Z_k\frac{\partial
  Z_m}{\partial \zeta_q}
\end{eqnarray}
pour tous $k=1,2,\ldots,n-p+1$,  $q=n-p+2,\ldots,n$ et
$m=1,\ldots,n$. En effet, par exemple pour $q=n$, $k=1$, 
si on fixe $a$, $\zeta_s$ pour tout $s=n-p+2,\ldots, n-1$ et
$\eta^r_s$ pour tout $(r,s)\not=(1,n)$, alors tout $Z_i$ s'\'ecrit
en fonction affine de $Z_n$. Il suffit donc de prouver l'\'egalit\'e
pr\'ec\'edente pour $m=n$. On \'ecrit $Z_1=\alpha Z_n+\beta$ et
$\eta^2_n Z_2+\cdots + \eta^{n-p+1}_n Z_{n-p+1}=\omega
Z_n+\theta$. On a 
$$Z_n=\zeta_n+\eta^1_n(\alpha Z_n+\beta)+\omega
Z_n+\theta.$$
D'o\`u
$$Z_n=\frac{\zeta_n+\eta^1_n\beta+\theta}{1-\eta^1\alpha-\omega}
\mbox{ et }Z_1=
\frac{\alpha\zeta_n+\alpha\theta+\beta-\beta\omega}{1-\eta^1\alpha-\omega}$$
$$\frac{\partial Z_n}{\partial
  \eta^1_n}=\frac{\alpha\zeta_n+\alpha\theta+\beta-\beta\omega}
{(1-\eta^1\alpha-\omega)^2}=Z_1\frac{\partial Z_n}{\partial \zeta_n}.$$
Les \'egalit\'es (3) se trouvent \'egalement dans
\cite[lemme de Darboux]{Henkin}. Elles nous donnent
aussi
\begin{eqnarray}
\frac{\partial \overline Z_m}{\partial \overline \eta^k_q} & = &
  \overline Z_k\frac{\partial
  \overline Z_m}{\partial \overline\zeta_q}
\end{eqnarray}
Dans les formules suivantes, les $z_k$ peuvent \^etre confondus avec
 leurs transformations
  d'Abel-Radon $Z_k$.  
D'apr\`es (3), (4), nous avons pour tout polyn\^ome $P$ en $z'$ et pour tout
$s,k=1,2,\ldots, n-p+1$
\begin{eqnarray}
\frac{\partial R(gPdz_s)}{\partial \eta^k_m} & = & 
\sum_{q=n-p+2}^{n}\left\{\int_\Upsilon \frac{\partial g}{\partial
    z_q}\frac{\partial Z_q}{\partial\eta^k_m}PdZ_s + \int_\Upsilon
  g\frac{\partial P}{\partial z_q}\frac{\partial Z_q}{\partial\eta^k_m}
  dZ_s\right\} + \nonumber\\
& & +\int_\Upsilon gPd\frac{\partial Z_s}{\partial \eta^k_m\nonumber}\\
& = & \sum_{q=n-p+2}^{n}\left\{\int_\Upsilon \frac{\partial g}{\partial
    z_q}\frac{\partial Z_q}{\partial\zeta_m}PZ_kdZ_s + \int_\Upsilon
  g\frac{\partial P}{\partial z_q}\frac{\partial Z_q}{\partial\zeta_m}
  Z_kdZ_s \right\} +\nonumber\\
& &  +\int_\Upsilon gPdZ_k\frac{\partial Z_s}{\partial \zeta_m}\nonumber\\
& = & \sum_{q=n-p+2}^{n}\left\{\int_\Upsilon \frac{\partial g}{\partial
    z_q}\frac{\partial Z_q}{\partial\zeta_m}PZ_kdZ_s + \int_\Upsilon
  g\frac{\partial P}{\partial z_q}\frac{\partial Z_q}{\partial\zeta_m}
  Z_kdZ_s\right\}+\nonumber\\
& & +\int_\Upsilon gPZ_k d\frac{\partial Z_s}{\partial \zeta_m}
 +\int_\Upsilon gP\frac{\partial Z_s}{\partial \zeta_m}dZ_k\nonumber\\
& = & \sum_{q=n-p+2}^{n}\left\{\int_\Upsilon \frac{\partial g}{\partial
    z_q}\frac{\partial Z_q}{\partial\zeta_m}PZ_kdZ_s + \int_\Upsilon
  g\frac{\partial P}{\partial z_q}\frac{\partial Z_q}{\partial\zeta_m}
  Z_kdZ_s\right\}+\nonumber\\
& & +\int_\Upsilon gPZ_k d\frac{\partial Z_s}{\partial \zeta_m}
 +\int_\Upsilon gP\frac{\partial Z_k}{\partial\zeta_m}dZ_s-\nonumber\\
& & -\int_\Upsilon gP\frac{\partial Z_k}{\partial\zeta_m}dZ_s
 +\int_\Upsilon gP\frac{\partial Z_s}{\partial \zeta_m}dZ_k\nonumber\\
& = & \frac{\partial R(gPz_kdz_s)}{\partial\zeta_m} 
-\int_\Upsilon gP\frac{\partial Z_k}{\partial\zeta_m}dZ_s
 +\int_\Upsilon gP\frac{\partial Z_s}{\partial \zeta_m}dZ_k\nonumber\\
& = & \frac{\partial R(gPz_kdz_s)}{\partial\zeta_m} 
+\int_\Upsilon gP\left\{ -\frac{\partial Z_k}{\partial\zeta_m}dZ_s
 +\frac{\partial Z_s}{\partial \zeta_m}dZ_k\right\}
\end{eqnarray} 
et 
\begin{eqnarray}
\frac{\partial R(gPdz_s)}{\partial \overline\zeta_m} & = &
\int_\Upsilon\sum_{q=n-p+2}^n\frac{\partial g}{\partial \overline
    z_q}\frac{\partial \overline Z_q}{\partial\overline \zeta_m} P
  dZ_s 
\end{eqnarray}
\begin{eqnarray}
\frac{\partial R(gPdz_s)}{\partial \overline\eta^k_m} & = &
\int_\Upsilon\sum_{q=n-p+2}^n\frac{\partial g}{\partial \overline
    z_q}\frac{\partial \overline Z_q}{\partial\overline \zeta_m} P
  \overline Z_k dZ_s 
\end{eqnarray}
En utilisant le changement des coordonn\'ees $(z',z'')\longmapsto
(z'-z''-\nu_0(1,z'))$, on se ram\`ene au cas o\`u $\nu_0=0$.\\
Soient $L_\lambda$ les vecteurs $\mathbb{C}$-ind\'epentdants de
$\mathbb{C}\otimes_{\mathbb{R}} H$ 
pour $\lambda=1,2,\ldots u$. Alors il existe des
matrices complexes 
$$b_\lambda=\left(
\begin{array}{c}
b_{n-p+2,\lambda}\\
\vdots\\
b_{n,\lambda}
\end{array}
\right)
\mbox{ }
c_\lambda= \left(
\begin{array}{cccc}
c^1_{n-p+2,\lambda} & c^2_{n-p+2,\lambda} & \cdots & c^{n-p+1}_{n-p+2,\lambda}\\
\vdots & \vdots & \ddots & \vdots\\
c^1_{n,\lambda} & c^2_{n,\lambda} & \cdots & c^{n-p+1}_{n,\lambda}\\
\end{array}
\right)
$$
$$
d_\lambda=\left(
\begin{array}{c}
d_{n-p+2,\lambda}\\
\vdots\\
d_{n,\lambda}
\end{array}
\right)
\mbox{ }
e_\lambda= \left(
\begin{array}{cccc}
e^1_{n-p+2,\lambda} & e^2_{n-p+2,\lambda} & \cdots & 
e^{n-p+1}_{n-p+2,\lambda}\\
\vdots & \vdots & \ddots & \vdots\\
e^1_{n,\lambda} & e^2_{n,\lambda} & \cdots & e^{n-p+1}_{n,\lambda}\\
\end{array}
\right)
$$
telles que 
\begin{eqnarray}
L_\lambda & = & \sum_{m=n-p+2}^n b_{m,\lambda}
\frac{\partial}{\partial\zeta_m} + 
\sum_{k=1}^{n-p+1}\sum_{m=n-p+2}^n
c^k_{m,\lambda}\frac{\partial}{\partial\eta^k_m}+ \nonumber\\ 
& & +\sum_{m=n-p+2}^n d_{m,\lambda}\frac{\partial}{\partial\overline\zeta_m} + 
\sum_{k=1}^{n-p+1}\sum_{m=n-p+2}^n
e^k_{m,\lambda}\frac{\partial}{\partial\overline\eta^k_m}.
\end{eqnarray}
Posons 
$$Q_{m,\lambda}=b_{m,\lambda}+\sum_{k=1}^{n-p+1}c^k_{m,\lambda}z_k$$
$$K_{m,\lambda}=d_{m,\lambda}+\sum_{k=1}^{n-p+1}e^k_{m,\lambda}
\overline z_k.$$  
Nous avons donc
\begin{eqnarray}
L_\lambda R(gPdz_s) & = & \sum_{m=n-p+2}^n \frac{\partial
R(gPQ_{m,\lambda}dz_s)}{\partial\zeta_m}+\nonumber\\
& & 
+\sum_{m=n-p+2}^n\sum_{k=1\atop k\not = s}^{n-p+1} 
c^k_{m,\lambda}\int_\Upsilon gP
\left\{ -\frac{\partial Z_k}{\partial\zeta_m}dZ_s
 +\frac{\partial Z_s}{\partial
   \zeta_m}dZ_k\right\}+\nonumber \\
& & 
+\sum_{m=n-p+2}^n \int_\Upsilon\sum_{q=n-p+2}^n 
\frac{\partial g}{\partial\overline z_q}\frac{\partial
    \overline Z_q}{\partial\overline\zeta_m}P K_{m,\lambda}dZ_s 
\end{eqnarray}
Supposons que la
condition 1 de la proposition 1 
n'est pas satisfaisante. On appelle $J'\subset J$
l'ensemble de tout $j'\in J$ tel que la condition 1 est valable pour
$z\in\gamma_{j'}$. Soit $j\in J$ tel que $\gamma_j$ n'est pas incluse dans
$bS_{j''}$ pour tout $j''\not\in J'\cup\{j\}$. Supposons que
$\gamma_j$ ne v\'erifie pas la condition 2.\\  
Posons 
$$Q=(Q_{n-p+2,\lambda},\ldots,Q_{n,\lambda})_{\lambda=1}^{u'}$$
$$Q'=(Q_{n-p+2,\lambda},\ldots,Q_{n,\lambda})_{\lambda=1}^{p-1}$$
$$Q''=(Q_{n-p+2,\lambda},\ldots,Q_{n,\lambda})_{\lambda=p}^{u'}$$
$$K=(K_{n-p+2,\lambda},\ldots,K_{n,\lambda})_{\lambda=1}^{u'}$$
$$K'=(K_{n-p+2,\lambda},\ldots,K_{n,\lambda})_{\lambda=1}^{p-1}$$
$$K''=(K_{n-p+2,\lambda},\ldots,K_{n,\lambda})_{\lambda=p}^{u'}$$
o\`u 
$u'=u-(p-1)(n-p)\geq 2p-1$.\\
On prend $s=1$. On choisit les $L_\lambda$ tels 
que $c^k_{m,\lambda}=0$ pour
tout $k\not =1$ et tout $\lambda=1,\ldots ,u'$. 
En plus, on peut choisir les $L_\lambda$
tels que $\rang Q'=p-1$ (c.-\`a-d. son d\'eterminant est un polyn\^ome
non nul).\\
On pose
$$T=(-\det Q'. Q''\times {Q'}^{-1},\det Q'. I_{u'-p+1})$$
o\`u $I_{u'-p+1}$ est la matrice d'identit\'e d'ordre $u'-p+1$.
Alors $T$ est de taille $(u'-p+1)\times
u'$ \`a coefficients dans $\P$. On a $T\times Q=0$.
\begin{definition} $H$ est un $u$-plan {\it g\'en\'erique} s'il existe 
un syst\`eme des coordonn\'ees $z'$, un choix des vecteurs $L_\lambda$
et une matrice $N$ form\'ee par
$p$ lignes de $T\times K$ tels que le d\'eterminant de la
matrice $(N,W)$ est non nul pour toute $p\times 1$-matrice non
nulle  $W$ \`a coefficients dans $\P$.
\end{definition}
\begin{lemme}
\begin{trivlist} \item[]{\rm 1.} 
L'ensemble des $u$-plans g\'en\'eriques passant
  par $0$ est un ouvert de Zariski, 
  non vide de la grassmannienne r\'eelle
  $\G_{\mathbb{R}}(u,\mathbb{R}^{2(p-1)(n-p+2)})$, qui param\`etre les
  $u$-plans passant par $0$.
\item[]{\rm 2.} Si $H$ est g\'en\'erique le d\'eterminant de 
  $(N,W)$ est non nul pour toute
  $p\times 1$-matrice non nulle $W$ 
\`a  coefficients dans $\K_j$.
\end{trivlist} 
\end{lemme}
\begin{preuve}
1. Soit $H$ un $u$-plan non g\'en\'erique. Les coefficients de la matrice
$T\times K$ sont des polyn\^omes de d\'egr\'es $1$ de $\overline z'$
\`a coefficients dans $\P$. 
Soit $N$ une matrice form\'ee par $p$-lignes de $T\times K$.\\
Posons $N^{(i)}$ la sous-matrice de $N$ obtenue par suppression de la
$i$-i\`eme ligne. Le d\'eterminant de $N^{(m)}$ s'\'ecrit sous la
forme 
$$\det N^{(i)}=\sum_{m=(m_1,\ldots,m_{n-p+1})\in\mathbb{N}^{n-p+1}
  \atop m_1+\ldots+ m_{n-p+1}\leq p-1}\sigma_{m,i}.\overline
  z_1^{m_1}\ldots  \overline z_{n-p+1}^{m_{n-p+1}}$$
o\`u $\sigma_{m,i}\in \P$.\\
Alors la matrice $\Sigma$ form\'ee par les $\sigma_{m,i}$ est de rang $\leq
  p-1$. Ceci implique le d\'eterminant de tout $p\times p$-sous-matrice de
  cette matrice $\Sigma$ est un polyn\^ome nul. Le $u$-plan $H$ est non
  g\'en\'erique si cette condition est v\'erifi\'ee pour toute
  $N$. Par cons\'equent, l'ensemble des
  $u$-plans g\'en\'eriques est un ouvert de Zariski.\\
Cet ouvert est non vide car, par exemple, si $H$ contient le plan
  engendr\'e par les vecteurs
\begin{enumerate}
\item $L_1=\displaystyle\frac{\partial}{\partial\zeta_{n-p+2}}$
\item $L_p=\displaystyle\frac{\partial}{\partial\overline\zeta_{n-p+2}}$
\item $L_{m-n+p-1}=\displaystyle\frac{\partial}{\partial\zeta_m}+
\frac{\partial}{\partial\overline
  \zeta_m} \mbox{ pour } m= n-p+3,\ldots,n$ 
\item $L_{m-n+2p-1}=\displaystyle\frac{\partial}{\partial\eta_m} + 
\frac{\partial}{\partial\overline \eta_m}  
+\frac{\partial}{\partial\eta_{m+1}} -
\frac{\partial}{\partial\overline \eta_{m+1}} \mbox{ pour }
  m=n-p+2,\ldots, n-1$
\item $L_{2p-1}=\displaystyle\frac{\partial}{\partial\eta_n} + 
\frac{\partial}{\partial\overline \eta_n}$
\end{enumerate}
Pour ces vecteurs on obtient 
$$N=\left(\begin{array}{cccccc}
 1 & 0 & 0 & \ldots  & 0 & 0 \\
\overline z_1 & -\overline z_1-z_1 & 0 & \ldots  & 0 & 0 \\
0 & \overline z_1 - z_1 & -\overline z_1-z_1 & \ldots   & 0 & 0 \\
\vdots & \vdots & \vdots & \ddots & \vdots & \vdots  \\
0 & 0 & 0 & \ldots & -\overline z_1-z_1 & 0 \\
0 & 0 & 0 & \ldots &  \overline z_1 -z_1 & -\overline z_1-z_1 \\
0 & 0 & 0 & \ldots &  0 & \overline z_1 -z_1
\end{array} \right)$$
Soit $W=^t(W_1,\ldots W_p)$ une $p\times 1$-matrice \`a coefficients
dans $\P$. Supposons que $\det (N,W)=0$. Alors 
$$
(-1)^pW_1 \overline z_1(\overline z_1 -z_1)^{p-2}+ \sum_{i=2}^p (-1)^{p-i+1}
W_i(\overline z_1 -z_1)^{p-i}(-\overline z_1-z_1)^{i-2}=0.$$
Comme $W_i\in \P$, ils sont tous nuls.\\
2. Comme dans 1., $W$ est d\'etermin\'ee comme un vecteur propre
d'une certaine matrice \`a coefficients dans $\P$. D'autre part, $\P$
s'injecte dans l'anneau int\`egre $\K_j$. Donc cette matrice n'a pas
de vecteur propre non nul \`a coefficients 
dans $\P$ si et seulement si elle n'a pas
de vecteur propre non nul \`a coefficients dans $\K_j$. D'o\`u 2.
\end{preuve}
Supposons que $H$ est g\'en\'erique. Supposons que la matrice $N$
form\'ee par les $p$ premi\`eres lignes de $T\times K$ v\'erifie la
d\'efinition 2.
Rempla\c cant $P$ par $PT_{i,\lambda}$ dans (9) et prenant la somme
en $\lambda=1,\ldots,2p-1$, nous obtenons
\begin{eqnarray}
\sum_{m=n-p+2}^n\int_\Upsilon\sum_{q=n-p+2}^n \frac{\partial
  g}{\partial\overline z_q}\frac{\partial\overline
  Z_q}{\partial\overline \zeta_m}P(T\times K)_{i,m}dz_1 & = & 0
\end{eqnarray}
car les d\'eriv\'ees $L_\lambda R(gPdz_1)=0$ par hypoth\`ese. Ceci montre que
\begin{eqnarray}
W_i:=\sum_{m=n-p+2}^n\left(\sum_{q=n-p+2}^n \frac{\partial
  g}{\partial\overline z_q}\frac{\partial\overline
  Z_q}{\partial\overline \zeta_m}\right)(T\times K)_{i,m}
\end{eqnarray}
appartient \`a $\K_j$ pour tout $i=1,\ldots, u'$ 
car son produit avec $dz_1$ est une mesure
  orthogonale. Comme $z''$ est un syst\`eme des coordonn\'ees de
$\mathbb{C}^{p-1}_a$ pour tout $a\in\Upsilon$, le d\'eterminant de la
matrice 
$$\left(\frac{\partial \overline Z_q}{\partial\zeta_m}
\right)_{q,m={n-p+2}}^n$$
est non nul. Par cons\'equent, $W\not =0$ car, d'apr\`es la
d\'efinition 2, $\rang N=p-1$.
D'autre part, le d\'eterminant de $(N,W)$
  est un polyn\^ome non nul en $\overline z'$ \`a coefficient dans
  $\K_j$. Par le th\'eor\`eme d'unicit\'e, on sait que 
  toute fonction de $\K_j$
  s'annullant sur un sous-ensemble de longueur positive de $\gamma_j$
  est nulle. Par cons\'equent, $\det (N,W)$ est de d\'egr\'e $\geq 1$.
  Ce polyn\^ome doit s'annuller sur $\gamma_j$ pour que
  (11) admette une solution non nulle pour $z\in\gamma_j$. C'est une
  contradiction.\\
Pour le cas o\`u $B$ n'est pas une r\'eunion de $(p-1)$-plans, on
consid\`ere $B'$ la r\'eunion de $(p-1)$-plans complexes tangentes \`a
$B$ en $z\in\gamma$. Les calculs pr\'ec\'edents sont encore valides
car ils ne concernent que des d\'eriv\'ees d'ordre 1, de plus, $B$
et $B'$ se co\"{\i}ncident en $\gamma$ \`a l'ordre 1.
\end{preuve}
\begin{definition}
Un sous-ensemble $V\subset
\G_D(n-p+2,n+1)$ est appel\'e {\it $(D,u)$-gros} si la r\'eunion
 $\bigcup_{\nu\in V'}\mathbb{C}^{n-p+1}_\nu\cap bD$ est dense dans
  $bD$, o\`u $V'$ est l'ensemble des \'el\'ements $\nu\in
  V$  tels que
  le c\^one tangent de $V$ en $\nu$ contienne au moins $u$ vecteurs
  r\'eellement ind\'ependants.
\item 
  Cet ensemble $V$ est appel\'e {\it
    $(D,u)$-g\'en\'erique} si il est $(D,u)$-gros avec $u\geq
  (p-1)(n-p+2)+1$ et 
$\overline{\bigcup_{\nu\in V''}\mathbb{C}^{n-p+1}_\nu\cap bD}$ est
dense dans 
  $bD$, o\`u $V''$ est l'ensemble des \'el\'ements $\nu\in
  V'$ v\'erifiant les conditions suivantes:
\begin{enumerate}
\item Le plan r\'eel engendr\'e par les vecteurs tangents \`a $V$ en
  $\nu$ est g\'en\'erique. 
\item $bD\cap \mathbb{C}_\nu^{n-p+1}$ est g\'en\'erique
  dans $\mathbb{C}_\nu^{n-p+1}$.
\end{enumerate}
\end{definition}
\begin{theoreme} Soient $Y$ un compact $(n-p+1)$-lin\'eairement convexe de
$\mathbb{CP}^n$ et $D$ une vari\'et\'e complexe (\'eventuellement
r\'eductible et singuli\`ere) de dimension $p\geq 2$,  \`a bord $\C^2$ dans
$\mathbb{C}^n\setminus Y$ et born\'ee dans $\mathbb{C}^n$. 
Soient $f$ une fonction $\C^1$ sur $bD$ et $V$ un ensemble
$(D,u)$-g\'en\'erique de $\G_D(n-p+2,n+1)$ avec $u\geq (p-1)(n-p+2)+1$. 
Supposons que pour tout $\nu\in V$, la fonction $f$ se prolonge
contin\^ument en une fonction holomorphe  de
$D\cap\mathbb{C}^{n-p+1}_\nu$. 
Alors $f$ est CR. Par cons\'equent, $f$ se prolonge contin\^ument en
une fonction holomorphe dans $D\setminus\Sing D$, localement born\'ee
dans $\overline D\setminus Y$. 
\end{theoreme}
\begin{remarque} Si dans la d\'efinition des ensembles
  $(D,u)$-g\'en\'eriques on supprime la condition 1, le th\'eor\`eme 2
  n'est plus valable ({\it voir} l'exemple 3). Si $u=2v$ et $V$ 
  une vari\'et\'e complexe de dimension $v$ de $\G(n-p+2,n+1)$, alors
  il suffit que $v\geq p$ \`a la condition que cet ensemble $V$
  soit g\'en\'erique dans un sens analogue que celui de la d\'efinition
  3.
\end{remarque}
\begin{preuve} D'apr\`es la proposition 1, la fonction $f$ v\'erifie
  la condition (1) sur un ensemble dense de $bD$. Par continuit\'e,
  (1) est satisfaisante partout sur $bD$. La fonction $f$ est donc
  CR. D'apr\`es le th\'eor\`eme de Dolbeault-Henkin g\'en\'eralis\'e 
  \cite[th\'er\`eme 1]{Dinh4} appliqu\'e pour le graphe
  de $f$, la propri\'et\'e
  ``1-extension'' de $f$ pour une grande famille de $(n-p+1)$-plans implique
  que le graphe de $f$ borde une vari\'et\'e complexe qui contient les
  graphes des prolongements de $f$ sur les tranches
  $\mathbb{C}^{n-p+1}_\nu$. Cette vari\'et\'e est le graphe d'un
  prolongement de $f$ sur $D\setminus\Sing D$ en fonction holomorphe
  localement born\'ee dans $\overline D\setminus Y$.  
\end{preuve}
\begin{lemme} Tout hyperplan r\'eel $H$ de
  $\mathbb{C}^{(p-1)(n-p+2)}$ passant par $0$ (c.-\`a-d. pour 
$u=2(p-1)(n-p+2)-1$) est g\'en\'erique.
\end{lemme}
\begin{preuve} On remarque que $H$ coupe tout sous-espace de
  $\mathbb{C}^{(p-1)(n-p+2)}$ en un sous-espace de codimension r\'eel
  $\leq 1$. Par cons\'equent, pour un syst\`eme des coordonn\'ees
  g\'en\'erique, on peut choisir les vecteurs $L_\lambda$ ind\'ependants
  v\'erifiant
\begin{enumerate}
\item Pour tout $\lambda\leq 4p-5$, pour tout $k\not=1$, pour tout $m$,
  $c^k_{m,\lambda}=e^k_{m,\lambda}=0$.
\item Pour tout $\lambda\leq p-2$, pour tout $m$, $c^1_{m,\lambda}=0$.
\item Pour tout $\lambda\leq 2p-3$, pour tout $m$ et tout $k$,
  $d_\lambda=e^1_{m,\lambda}=0$.
\item Pour tout $p-1\leq\lambda\leq 2p-3$, pour tout
  $m\not=\lambda+n-2p$, $b_{m,\lambda}=c^1_{m,\lambda}=0$.
\item  Pour tout $p-1\leq\lambda\leq 2p-3$,
  $c^1_{\lambda+n-2p,\lambda}=1$.
\item Pour $2p-2\leq\lambda\leq 4p-6$, pour tout $m$, $k$,
  $b_{m,\lambda}=c^1_{m,\lambda}=0$, $d_{m,\lambda}=\overline
  b_{m,\lambda-2p+3}$ et $e^k_{m,\lambda}=\overline
  e^k_{m,\lambda-2p+3}$.
\item Pour tout $m$, $c^1_{m,4p-5}=e^1_{m,4p-5}=0$.
\item Le d\'eterminant de $Q'$ est un polyn\^ome non nul $F$ de
  d\'egr\'e $\leq 1$ en $z_1$.
\item $b_{m,4p-5}=\overline d_{m,4p-5}\not=0$ pour tout $m$.
\end{enumerate}
Les $(4p-6)$ premiers vecteurs forment une base de l'espace
$\mathbb{C}\otimes_{\mathbb{R}}H'$, o\`u $H'$ est 
le sous-espace complexe maximal de $H\cap \{\eta^k_m=0\mbox{ pour tout
  } k\not=1\}$.\\
Les $2p-3$ premi\`eres lignes de $T\times K$ et de $K''$
sont \'egales, dont les $p-2$ premi\`eres lignes sont nulles. La
derni\`ere  ligne de $T\times K$ est $F$ fois celle de
$K''$. Supposons que $H$ n'est pas g\'en\'erique, alors il n'existe
pas de matrice $N$ v\'erifie la d\'efinition 2 pour toute combinaison
lin\'eaire \`a coefficients dans $\mathbb{C}$ des vecteurs $L_\lambda$
pour $\lambda=2p-2,\ldots,4p-5$. On sait que les lignes num\'eros
$\lambda=2p-2,\ldots ,3p-4$ et $\lambda=4p-5$ de $K$ forment une base de
$\mathbb{C}^{p-1}$. En rempla\c cant  $L_\lambda$
pour $\lambda=2p-2,\ldots,4p-6$ par leurs combinaisons convenables, on
peut trouver la matrice $N$ form\'ee par la derni\`ere ligne de
$T\times K$ (not\'ee $FN_1$) et par les lignes de la forme
$$N_1+(0,\ldots,0,k+\overline z_1,0,\ldots,0)$$
pour $k=1,\ldots, p-1$ et $k+\overline z_1$ est \`a la $k$-i\`eme
position. Comme la matrice $N$ ne v\'erifie pas la d\'efinition 2, la
matrice $N'$ form\'ee par 
$$N_1=(d_{n-p+2,4p-5},\ldots,d_{n,4p-5})$$ 
et par les autres lignes de $N$, ne v\'erifie pas la d\'efinition 2
non plus. Par cons\'equent, la matrice $N''$ form\'ee par $N_1$ et les
lignes de la forme
$$(0,\ldots,0,k+\overline z_1,0,\ldots,0)$$
ne v\'erifie pas la d\'efinition 2. Posons
$A=\Pi_{k=1}^{p-1}(k+\overline z_1)$. Alors 
il existe $W$ non nulle \`a coefficients dans $\P$ telle que
$$W_1A\pm W_2 d_{n-p+2,4p-5}\frac{A}{1+\overline z_1}\pm\ldots\pm W_p
d_{n,4p-5}\frac{A}{p-1+\overline z_1}=0$$
Ceci implique $W=0$. C'est une contradiction.
\end{preuve}
\begin{corollaire}
Le compact $Y$, la vari\'et\'e $D$ et la fonction 
$f$ sont d\'efinis comme dans le th\'eor\`eme 2. 
Soit $V\subset\G_D(n-p+2,n+1)$ tel que
$bD\cap\bigcup_{V_1}\mathbb{C}^{n-p+1}_\nu$ soit dense dans $bD$, o\`u
$V_1\subset V$ est l'ensemble de $\nu$ v\'erifiant les conditions
suivantes:
\begin{enumerate}
\item Le c\^one tangent de $V$
en $\nu$ contient $2(p-1)(n-p+2)-1$ vecteurs ind\'ependants.
\item L'intersection $\mathbb{C}^{n-p+1}_\nu\cap bD$ est une courbe
  g\'en\'erique dans $\mathbb{C}^{n-p+1}_\nu$.
\end{enumerate}
Supposons que pour tout
  $\nu\in V$, $f$ se prolonge contin\^ument en une fonction holomorphe 
  dans $\mathbb{C}^{n-p+1}_\nu\cap D$. Alors $f$ se prolonge contin\^ument en
  une fonction holomorphe 
  dans $D\setminus\Sing D$, localement
  born\'ee dans $\overline D\setminus Y$.
\end{corollaire}
\begin{lemme}
Soient $\gamma$ une courbe r\'eelle, ferm\'ee, irr\'eductible bordant un
domaine $\Omega\subset \mathbb{C}$ et $l$ un ouvert de $\gamma$. 
Supposons qu'il existe un voisinage $U$ de $l$ et 
deux fonctions holomorphes non identiquement nulles 
$a,b$  dans $U\cap \Omega$, continue jusqu'\`a tout point de $l$ et
v\'erifiant $a\overline z +b=0$ sur $l$. Alors
$l$ est presque r\'eelle analytique. Si $l=\gamma$ est r\'eelle
analytique et
$U\supset \Omega$, alors $\gamma$ est r\'eelle alg\'ebrique (voir
l'exemple 1).
\end{lemme}
\begin{preuve} Soit $K\subset l$ l'ensemble des z\'eros de $a$ sur
  $l$.
D'apr\`es le
  th\'eor\`eme d'unicit\'e $\H^1(K)=0$. 
  Soient $x\in l\setminus K$, 
  $G$ un domaine de Jordan \`a bord lisse de $\mathbb{C}$ dont le bord
  contient le segment r\'eel
  $[-1,1]\subset\mathbb{R}\subset\mathbb{C}$ et $\Psi$ une application
  bijective de $\overline G$ dans $\Omega$, holomorphe dans $G$
  telle que $\Psi([-1,1])\subset l$ et $\Psi(0)=x$. Comme
  $\overline z=-b/a$ on a $|z|^2=-bz/a$. Par cons\'equent, la fonction
  $|\Psi|^2$ s'\'etend m\'eromorphiquement dans $G$ (holomorphe
  sur les points pr\`es de $0$). D'apr\`es le principe de r\'eflexion,
  cette fonction s'\'etend en une fonction holomorphe au voisinage de
  $0$. Donc $l$ est analytique au voisinage de $x$. \\
Si $l=\gamma$ r\'eelle analytique, on prend $G$ le demi-plan complexe
sup\'erieur et $\Psi$ l'application bijective de $\overline
G\subset\mathbb{CP}^1$ dans $\overline\Omega$, holomorphe dans
$G$. Alors $|\Psi|^2$ se prolonge holomorphiquement au voisinage de
$bG$ et \'egale \`a $-a/b$ $\H^1$-presque partout sur $bG$. Par
cons\'equent, cette fonction se prolonge contin\^ument en fonction
m\'eromorphe dans $G$. D'apr\`es le principe de r\'eflexion, elle se
prolonge en une fonction m\'eromorphe dans $\mathbb{CP}^1$,
c.-\`a-d. en une fonction rationnelle. On a $|\Psi|^2=P/Q$, o\`u $P,Q$
sont des polyn\^omes \`a coefficients r\'eels. 
On peut supposer que
$0\not\in \overline\Omega$. Alors il existe une fonction $\Phi$
d\'efinie sur $bG$ \`a valeur dans $\{|z|=1\}$ 
telle que sur $bG$ on ait $\Psi^2=\Phi P/Q$. La fonction
$\Phi$ s'\'etend m\'eromorphiquement dans $G$. 
Par principe de r\'eflexion, elle s'\'etend
m\'eromorphiquement dans $\mathbb{CP}^1$. Ce prologement est donc une
fonction rationnelle. Par cons\'equent, $b\Omega$ est r\'eelle alg\'ebrique. 
\end{preuve}
{\it Preuve du th\'eor\`eme 1.} Utilisant les $L_\lambda$ dans le
lemme 2 pour un $\nu\in V_1$ fix\'e, 
supposons que la condition 1 de la proposition 1 n'est pas
satisfaisante. Alors la condition 2 de cette proposition montre que
pour un syst\`eme des coordonn\'ees g\'en\'erique, il existe un
polyn\^ome en $\overline z_1$ \`a coefficients dans $\K_j$ s'annullant
sur $\gamma_j$. D'apr\`es le lemme 3, la projection $\Pi(\gamma_j)$
est presque r\'eelle analytique, o\`u $\Pi(z)=z_1$. Ceci est valable
pour un syst\`eme des coordonn\'ees g\'en\'erique. Par cons\'equent,
toute coordonn\'ee de $\gamma_j$ est presque r\'eelle analytique. Donc
$\gamma_j$ est presque r\'eelle analytique. C'est une
contradiction. Alors $f$ v\'erifie (1) pour tout $z\in bD\cap
\mathbb{C}^{n-p+1}_\nu$ et pour tout $\nu\in V_1$. Par continuit\'e, (1)
est vraie partout sur $bD$. La fonction $f$ est donc CR. D'apr\`es le
th\'eor\`eme de Dolbeault-Henkin g\'en\'eralis\'e
\cite{Dinh4}, $f$ s'\'etend holomorphiquement dans $D\setminus\Sing D$
en fonction localement born\'ee dans $\overline D\setminus Y$.
\\
$\square$
\begin{theoreme}
Soient 
$D$ un domaine de $\mathbb{C}^n$, born\'e, \`a bord $\C^2$ 
et $f$ une fonction $\C^1$ sur $bD$. Supposons que 
\begin{eqnarray}
\int_{bD\cap \mathbb{C}_{\nu}} 
f(z_1,\zeta_2+\eta_2 z_1\ldots,\zeta_n+\eta_n z_1)z_1^kdz_1 & = & 0
\end{eqnarray} 
pour $\H^{4n-4}$-presque toute $\nu=(\zeta,\eta)\in\mathbb{C}^{2n-2}$ 
et pour $k\in
\mathbb{N}$ fix\'e. Alors $f$ est CR et si $bD$ est connexe, 
$f$ se prolonge
holomorphiquement dans $D$.
\end{theoreme}
\begin{remarque} La condition (12) a \'et\'e introduite par  
Kytmanov et Myslivets dans 
\cite{KytmanovMyslivets}, o\`u ils ont g\'en\'eralis\'e des
r\'esultats de Globevnik et de Stout.
\end{remarque} 
\begin{preuve} Il suffit de consid\'erer
  $n=2$. On appelle $T$ (resp. $T'$) 
  l'ensemble de $\nu\in\mathbb{C}^2$ tel que
  $\mathbb{C}_\nu$ appartienne au plan tangent de $bD$ en certain
  point (resp. l'intersection $\mathbb{C}_\nu\cap bD$ n'est pas $\C^2$
  par morceaux). Alors l'ensemble $T$ (resp. $T'$) 
  est de mesure $\H^3$ (resp. $\H^2$) localement finie. 
  D'apr\`es (5), on a
$$\frac{\partial R(f z_1^{k-1}dz_1)}{\partial\eta_2}=\frac{\partial
  R(f z_1^k dz_1)}{\partial\zeta_2}=0$$
par hypoth\`ese. Donc sur la droite $\{\zeta_2=c\}$ 
la fonction $R(f z_1^{k-1}dz_1)$ est continue en dehors de $T'$ et 
antiholomorphe en dehors de $T$, o\`u $c\in\mathbb{C}$ une constante. 
Cette fonction est nulle pour $\eta_2$
  assez grand. Elle est donc nulle pour $\H^2$-presque tout
  $c$. D'apr\`es une r\'ecurrence en $k$, on conclure que $f$ v\'erifie la
  condition de Morera. De m\^eme mani\`ere, en consid\'erant les droites
  $\{\eta_2=c\}$, nous pouvons 
  d\'emontrer que $f$ v\'erifie la condition des moments. 
D'apr\`es le th\'eor\`eme de Stout, $f$ est CR et si $bD$ est connexe,
  la fonction $f$ se prolonge holomorphiquement sur $D$
  \cite{Stout1}.       
\end{preuve}
\begin{exemple} Soient $n=p=2$ et $\G^*_Y(2,3)$ un voisinage de
  $(0,0)$. Soient $P$ un polyn\^ome d'une variable complexe de
  d\'egr\'e $k$ \`a coefficients r\'eels avec
  $P(0)=0$, $P'(0)\not =0$ et $0<r<1$ tels que la restriction de $P$
  sur $\{|t|<r\}$ soit injective. Posons 
  $D=\mathbb{C}^2\setminus Y\cap\{z_1\in P(\{|t|=r\})\}$ et
  $f(z)=z_1^k \overline
  z_2$ une fonction d\'efinie sur $bD$. Alors pour tout
  $\nu=(\zeta,\eta)\in\G^*_Y(2,3)$ la droite $\mathbb{C}_\nu$ est
  d\'efinie par l'\'equation $\{z_2=\zeta+\eta z_1\}$. Sur
  $bD\cap\mathbb{C}_\nu$ on a $f(z)=z_1^k\overline z_2=[P(t)]^k(\overline\zeta
  + \overline\eta \overline{P(t)})=\overline \zeta [P(t)]^k+\overline
  \eta [P(t)]^kP(\overline t)= \overline \zeta [P(t)]^k+\overline
  \eta [P(t)]^kP(r^2/t)$ pour $z_1=P(t)$. 
  Cette fonction en $t$ s'\'etend holomorphiquement dans
  $\{|t|<r\}$. Par cons\'equent,
  $f$ s'\'etend holomorphiquement dans 
  $D\cap\mathbb{C}_\nu$ car dans $\{|t|<r\}$ l'application 
  $P$ est injective.  Mais la fonction $f$ n'est pas CR.
\end{exemple}
\begin{exemple} Pour $p=n=2$; $\nu_0=(0,0)$ et 
  $B\cap\mathbb{C}_{\nu_0}=\{z_1=\Psi(\e^{i\theta}),z_2=0\}$ o\`u
  $\Psi$ une fonction $\C^k$ d\'efinie sur le disque unit\'e ferm\'e
  $\overline U$
  avec
$$\Psi(t)=t+\epsilon (t-1)^{2k+1}\e^{\frac{t+1}{t-1}}.$$
(pour $\epsilon$ r\'eel assez petit et $k\geq 1$, 
cette application $\Psi$ est injective
sur $\overline U$). Sur   $B\cap\mathbb{C}_{\nu_0}$,
on a pour $t=\e^{i\theta}$
$$\overline z_1=\overline t +\epsilon(\overline t
-1)^{2k+1}\e^{\frac{\overline t+1}{\overline t-1}}=\frac{1}{t}+\epsilon
  \frac{(1-t)^{2k+1}}{t^{2k+1}} e^{\frac{1+t}{1-t}}$$
d'o\`u
$$\overline
z_1=\frac{1}{\Psi^{-1}(z_1)}-\epsilon^2
\frac{(1-\Psi^{-1}(z_1))^{4k+2}}{(\Psi^{-1}(z_1))^{2k+1}}
\frac{1}{z_1-\Psi^{-1}(z_1)}.$$
La restriction de $\overline z_1$ sur la courbe $B\cap\mathbb{C}_{\nu_0}$
s'\'etend m\'eromorphiquement dans le domaine de $\mathbb{C}$, born\'e
par cette courbe. Cette courbe est localement 
analytique sauf au point $z_1=1$.\\
Nous pouvons maintenant construire un exemple comme celui
pr\'ec\'edent. Soient $p=n=2$, $\G^*_Y(2,3)$ un voisinage de $(0,0)$,
  $D=\mathbb{C}^2\setminus Y\cap\Psi(U)\times \mathbb{C}$ et 
$$f(z)= (\Psi^{-1}(z_1))^{2k+1}(z_1-\Psi^{-1}(z_1))\overline z_2$$ 
une fonction d\'efinie sur $bD$. Pour
$\nu=(\zeta,\eta)\in\G^*_Y(2,3)$, on a sur $bD\cap\mathbb{C}_\nu$
$$f(z)=
(\Psi^{-1}(z_1))^{2k+1}(z_1-\Psi^{-1}(z_1))
(\overline\zeta+\overline \eta \overline z_1).$$
Cette fonction s'\'etend holomorphiquement sur
$D\cap\mathbb{C}_\nu$. Mais la fonction $f$ n'est pas CR.
\end{exemple}
\begin{exemple} Soient $n=p\geq 3$, $D$ un domaine born\'e dans
  $\mathbb{C}^n$, convexe, \`a bord $\C^2$ et $f(z)=\overline z_n$. 
  Soient $V=\{\eta^1_n=0\}
  \subset\G(2,n+1)$
  une sous-vari\'et\'e complexe de dimension $2n-3$. Alors $V$ est
  $(D,4n-6)$-gros. On peut choisir $D$ tel que $V$ v\'erifie la
  deuxi\`eme condition de la d\'efinition des ensembles 
  $(D,4n-6)$-g\'en\'eriques (d\'efinition 3). Pour tout
  $\nu\in V$, la fonction $f$ est constante sur
  $bD\cap\mathbb{C}_\nu\subset \{z_n=\zeta_n\}$. Elle s'\'etend donc
  holomorphiquement dans $bD\cap\mathbb{C}_\nu$. Mais cette fonction $f$
  n'est pas CR.
\end{exemple}
\section{Probl\`eme du bord}
Soit $\Gamma$ un courant rectifiable de
dimension $2p-1$ d'une vari\'et\'e complexe  $X=\mathbb{C}^n$, 
$\mathbb{CP}^n$, $\mathbb{C}^n\setminus Y$ ou $\mathbb{CP}^n\setminus
Y$. Le courant $\Gamma$ est appel\'e {\it maximalement complexe} si
pour toute $(s,2p-1-s)$-forme $\psi$ de classe $\C^\infty$, \`a
support compact dans $X$ 
et pour tout $s\not =p,p-1$, on a $(\Gamma,\psi)=0$. \\
{\it Une $p$-cha\^{\i}ne holomorphe de
$X\setminus\Supp\Gamma$} est une combinaison lin\'eaire localement finie 
\`a coefficients entiers de
sous-vari\'et\'es complexes de dimension pure $p$ de $X
\setminus\Supp \Gamma$. Si
une $p$-cha\^{\i}ne holomorphe $T$ est de mesure ${\cal H}^{2p}$,
compt\'ee avec la valeur absolue des coefficients, localement finie dans
$X$, 
elle d\'efinit dans $X$ un courant d'int\'egration de bidimension
$(p,p)$ et de masse localement finie.
\\
Si $X=\mathbb{C}^n$, pour  
${\cal H}^{(n-p+2)(p-1)}$-presque tout
$\nu\in\G(n-p+2,n+1)$ le courant d'intersection 
$\Gamma\cap \mathbb{P}^{n-p+1}_{\nu}$ existe et il 
est rectifiable \cite{Federer}.
\begin{theoreme} Soient $Y$ un compact de $\mathbb{C}^n$, 
$(n-p+1)$-lin\'eairement
  convexe dans $\mathbb{CP}^n$, $\Gamma$ un courant ferm\'e, de dimension
  $2p-1$, \`a
support ($\Supp\Gamma$) g\'eom\'etriquement $(2p-1)$-rectifiable 
dans $\mathbb{C}^{n}\setminus Y$ et born\'e dans $\mathbb{C}^n$, avec 
$p\geq 2$. Supposons que pour ${\cal H}^{(n-p+2)(p-1)}$-presque tout
$\nu\in\G_Y(n-p+2,n+1)$, $\Gamma\cap\mathbb{C}^{n-p+1}_\nu$ existe et  
v\'erifie la
condition de Morera ($(\Gamma\cap\mathbb{C}^{n-p+1}_\nu,z_idz_j)=0$
  pour tous $i,j$). Alors il existe une $p$-cha\^{\i}ne holomorphe
  $T$  de $\mathbb{C}^{n}\setminus\Supp\Gamma\cup Y$, de masse
  localement finie dans $\mathbb{C}^n\setminus Y$ telle que
$d[T]=\Gamma$ au sens des courants dans $\mathbb{C}^n\setminus Y$.
\end{theoreme}
\begin{remarque} Ce th\'eor\`eme g\'en\'eralise le th\'eor\`eme de
  Harvey-Lawson \cite{HarveyLawson} 
et il donne la r\'eponse \`a un probl\`eme de
  Dolbeault-Henkin \cite{DolbeaultHenkin}.\\ 
Si $Y=\emptyset$, il suffit de consid\'erer une
  famille de $(n-p+1)$-plans dont les directions appartiennent \`a un
  ouvert de $\G(n-p+1,n)$ ({\it voir} la preuve du lemme 4).
Ce th\'eor\`eme n'est
plus valable si l'on remplace $\mathbb{C}^{n}$ par $\mathbb{CP}^{n}$.
\end{remarque}
\begin{exemple} (Henkin) 
Soit $\Gamma\subset\mathbb{C}^{3}\subset\mathbb{CP}^{3}$
d\'efinie par
$$
\Gamma=
\{ y_{2}=y_{3}=0,x_{1}^{2}+y_{1}^{2}+x_{2}^{2}+x_{3}^{2}=1
\}
$$
o\`u $z_{1}=x_{1}+iy_{1}$, $z_{2}=x_{2}+iy_{2}$, $z_{3}=x_{3}+iy_{3}$
sont les coordonn\'ees de $\mathbb{C}^{3}$. Consid\'erons l'hyperplan
$$H_{a,b,c}=\{z_{1}=az_{2}+bz_{3}+c\}$$ 
o\`u $a=a_{1}+ia_{2}$,
 $b=b_{1}+ib_{2}$ et $c=c_{1}+ic_{2}$. Posons $\Gamma_{a,b,c}=\Gamma\cap
H_{a,b,c}$. Alors
$\Gamma_{a,b,c}$ est une courbe r\'eelle ferm\'ee de la surface de Riemann
alg\'ebrique $S_{a,b,c}\subset H_{a,b,c}$ qui est d\'efinie par
$$S_{a,b,c}=\{(a_{1}z_{2}+b_{1}z_{3}+c_{1})^{2} +
(a_{2}z_{2}+b_{2}z_{3}+c_{2})^{2} +z_{2}^{2}+z_{3}^{2}=1\}\cap H_{a,b,c}.
$$
Comme $S_{a,b,c}$ est de genre $0$, la courbe $\Gamma\cap H_{a,b,c}$ borde
une surface de Riemann dans $\mathbb{CP}^{3}$. 
La vari\'et\'e $\Gamma$ ne peut pas \^etre le
bord d'une vari\'et\'e complexe car elle n'est pas maximalement
complexe.
\end{exemple}
\begin{preuve}
D'apr\`es le th\'eor\`eme de Dolbeault-Henkin g\'en\'eralis\'e 
\cite[th\'eor\`eme 3]{Dinh4}, il suffit de montrer que
$\Gamma$ est maximalement complexe. Dans cette d\'emonstration on
consid\`ere $Y=\emptyset$, pour le cas g\'en\'eral, il suffit 
d'\'etudier le probl\`eme au voisinage
d'un $(n-p+1)$-plan fix\'e; en choisissant un bon syst\`eme des
coordonn\'ees, la m\^eme d\'emonstration pour $Y=\emptyset$
s'adaptera.\\ 
Par la m\'ethode des projections, il suffit de consid\'erer le cas
$p=n-1$.
Soit $\Pi$ une projection de $\mathbb{C}^{n}$ dans
$\mathbb{C}^{n-1}$ telle que sa restriction 
en $\Supp\Gamma$ soit injective sauf
au dessus d'un sous-ensemble $K$ de mesure ${\cal H}^{2n-3}$ nulle de
$\mathbb{C}^{n-1}$ et que la projection de $\Supp\Gamma$ soit 
g\'eom\'etriquement $(2p-1)$-rectifiable. En plus, $\Pi^{-1}(x)$ coupe
$\Tan(\Supp\Gamma,z)$ transversalement pour tout
$\Pi(z)=x\in\Pi(\Supp\Gamma)\setminus K$.
Ceci est valable une ${\cal H}^{2n-2}$-presque
toute projection $\Pi$ \cite{Dinh2}. Sans perdre en
g\'en\'eralit\'e, on suppose que $\Pi(z)=(z_{1},\ldots, z_{n-1})$ et
et pour ${\cal H}^{2n-4}$-presque toute projection $\Phi$ de
$\mathbb{C}^{n-1}$ dans $\mathbb{C}^{n-2}$ et
pour ${\cal H}^{2n-4}$-presque tout
$\nu\in\mathbb{C}^{n-2}$ le courant $\Gamma\cap\mathbb{C}^{2}_{\nu}$
existe et v\'erifie la condition du th\'eor\`eme 1, o\`u
$\mathbb{C}^{2}_{\nu}=(\Phi\circ\Pi)^{-1}(\theta)$. Le support
$\Supp\Gamma$ est d\'efini $\H^{2n-3}$-presque partout comme le graphe
d'une fonction $f$ au dessus de $\Pi(\Supp\Gamma)$.
Le lemme suivant est prouv\'e gr\^ace \`a l'utilisation d'une 
id\'ee de Globevnik-Stout \cite{GlobevnikStout}:
\begin{lemme} Soit $\Pi_{*}\Gamma^{0,1}$ la composante de bidegr\'e
$(0,1)$ du courant $\Pi_{*}\Gamma$. Alors le courant
$f\Pi_{*}\Gamma^{0,1}$ est $\overline{\partial}$-ferm\'e.
\end{lemme}
\begin{preuve} D'apr\`es \cite{GlobevnikStout}, l'espace de fonctions
${\cal C}^{\infty}$ de $\mathbb{C}^{n-1}$ \`a valeurs complexes 
engendr\'ees par des fonctions, dont
les lignes de niveau sont des hyperplans complexes parall\`eles, est dense
dans l'espace des fonction ${\cal C}^{\infty}$. (Plus g\'en\'eralement,
ce sous-espace est d\'ej\`a dense 
si on  consid\`ere seulement les hyperplans
parall\`eles dont les directions sont param\'etr\'ees par un ouvert
non vide de $\G(n-2,n-1)$).\\
Sans perdre en
g\'en\'eralit\'e, il suffit de prouver que
$(\Pi_{*}\Gamma,f\overline{\partial}\alpha)=0$ pour toute\break
$(n-1,n-3)$-forme $\alpha=A(\zeta \tilde z)\chi$
pour un certain $\zeta$ g\'en\'erique non nul de $\mathbb{C}^{n-1}$, 
o\`u $\tilde z=(z_1,\ldots,z_{n-1})$, $\zeta \tilde z=\zeta_1
z_1+\ldots \zeta_{n-1}z_{n-1}$, 
$\chi=dz_{1}\wedge\ldots\wedge
dz_{n-1}\wedge d\overline{z}_{1}\wedge\ldots\wedge d\overline{z}_{n-3}$,
$\zeta \tilde 
z=\zeta_{1}z_{1}+\cdots +\zeta_{n-1}z_{n-1}$ et $A$ est une fonction
${\cal C}^{\infty}$.\\
On utilise un nouveau syst\`eme de coordonn\'ees $(w_{1},\ldots,w_{n-1})$
avec $w_{1}=\zeta z$. Alors il existe des constantes $a_{k,j}\in
\mathbb{C}$ telles que $\chi=\sum_{1\leq j< k\leq p}a_{j,k}\chi_{j,k}$, o\`u 
$$
\chi_{j,k}=dw_{1}\wedge\ldots\wedge dw_{n-1}\wedge
d\overline{w}_{1}\wedge\ldots\wedge\widehat{d\overline{w}_{j}}\wedge\ldots
\wedge\widehat{d\overline{w}_{k}}\wedge\ldots\wedge\ d\overline{w}_{n-1}.
$$
On a
$$
\overline{\partial}A\chi_{j,k}  = 0 
$$
si $j\geq 2$ et
\begin{eqnarray*}
(\Pi_{*}\Gamma,f\overline{\partial}A\chi_{1,k}) & = & 
(\Pi_{*}\Gamma,f\frac{{\partial}A}{\partial\overline{w}_{1}} 
d\overline{w}_{1}\wedge \chi_{1,k})\\
& = &
\displaystyle
\pm\int_{a\in\mathbb{C}^{n-2}}\frac{{\partial}A(a_{1})}{\partial
\overline{a}_{1}} 
(\Pi_{*}\Gamma\cap \Psi_{k}^{-1}(a),f dw_{k})da_{1}\wedge\ldots\wedge
 d\overline{a}_{n-2}\\
& = & 0
\end{eqnarray*}
car $$(\Pi_{*}\Gamma\cap \Psi_{k}^{-1}(a),fdw_{k})=
(\Gamma\cap \Psi_{k}^{-1}(a),z_n dw_k)=0$$
par hypoth\`ese, 
o\`u $\Psi_{k}(w):=(w_{1},\ldots,\widehat{w_{k}},\ldots,w_{n-1})$
et $\nu\in \G(3,n)$ est d\'efini par 
$\mathbb{C}^2_\nu=\Psi_{k}^{-1}(a)\subset\mathbb{C}^{n-1}$.
\\
Ces \'egalit\'es donnent le lemme.
\end{preuve}
D'apr\`es Harvey, si $\Gamma$ est une
vari\'et\'e $\C^1$, le lemme pr\'ec\'edent implique que le courant
$\Gamma$ est maximalement complexe \cite{Harvey}. Cette proposition
est encore valable dans notre cas: le plan tangent de $\Supp\Gamma$
est maximalement complexe en $\H^{2n-2}$-presque tout point de
$\Supp\Gamma\setminus\Pi^{-1}(K)$ 
o\`u la multiplicit\'e de $\Gamma$ est non nul. En  
appliquant le lemme pr\'ec\'edent pour les projections diff\'erentes, on
conclure que $\Gamma$ est maximalement complexe. D'apr\`es le
th\'eor\`eme de Harvey-Lawson g\'en\'eralis\'e \cite{Dinh2}, 
$\Gamma$ est donc le bord
d'une $p$-cha\^{\i}ne holomorphe de masse finie au sens des courants.
\end{preuve}
\begin{theoreme} Soient $Y$ un compact $(n-p+1)$-lin\'eairement
  convexe de $\mathbb{CP}^n$, $\Gamma$ une vari\'et\'e $\C^2$ de
  dimension $2p-1$, orient\'ee de $\mathbb{C}^n\setminus Y$ et
  born\'ee dans $\mathbb{C}^n$. Soit $V$ une vari\'et\'e r\'eelle de
  dimension $2(p-1)(n-p+2)-1$ immerg\'ee dans
  $\G_Y(n-p+2,n+1)$. Supposons que 
\begin{enumerate}
\item Pour tout $\nu\in V$, $\Gamma\cap\mathbb{C}^{n-p+1}_\nu$ est
  transversale et 
  borde une $1$-cha\^{\i}ne holomorphe, de masse finie au sens des courants.
\item $\Gamma\cap\bigcup_{\nu\in
  V_1}\mathbb{C}^{n-p+1}_\nu$ est dense dans $B$, o\`u $V_1$ est
  l'ensemble de $\nu\in V$ tel que aucun ouvert non vide de
  $\Gamma\cap \mathbb{C}^{n-p+1}_\nu$ n'est r\'eelle analytique.
\end{enumerate} 
Alors $\Gamma$ est le bord d'une $p$-cha\^{\i}ne holomorphe de masse
localement finie au sens
des courants dans $\mathbb{C}^n\setminus Y$.
\end{theoreme}
\begin{preuve} D'apr\`es le th\'eor\`eme de Dolbeault-Henkin
  g\'en\'eralis\'e \cite{Dinh4}, il suffit de montrer que
  $\Gamma$ est maximalement complexe.
On fixe $\nu_0\in V_1$, il suffit de prouver que $\Gamma$ est
maximalement complexe en tout point de
$\Gamma\cap\mathbb{C}^{n-p+1}_{\nu_0}$. Par la m\'ethode de
projection, 
on  ram\`ene le probl\`eme
 vers le cas o\`u $n=p+1$. Soit $\Pi$ une projection g\'en\'erique de
 $\mathbb{C}^n$ dans $\mathbb{C}^{n-1}$ v\'erifiant
 $\Pi(\mathbb{C}^{2}_{\nu_0})\simeq \mathbb{C}$. Posons
 $\Lambda=\Pi(\mathbb{C}^{2}_{\nu_0})$. Alors au
 voisinage de $\Lambda$, $\Pi(\Gamma)$ est $\C^2$ par morceaux et
 $\Gamma$ est d\'efinie comme le graphe d'une fonction born\'ee $f$ au dessus
 de $\Pi(\Gamma)$ (sauf sur les singularit\'es). La
 proposition 1 s'applique encore dans ce cas,  $f$ v\'erifie donc (1) 
en tout point r\'egulier de $\Pi(\Gamma)\cap\Lambda$.
Par cons\'equent, $\Gamma$ est
 maximalement complexe en tout point de 
$\Gamma\cap\mathbb{C}^{2}_{\nu_0}$.   
\end{preuve}

\begin{thebibliography}{11}
\bibitem{AgranovskiSemenov}
\textit{M.L. Agranovski, A.M. Semenov}, Boundary analogues of
Hartog's theorem, Sibirian. Math. J., \textbf{32} (1991), 168-170.
\bibitem{Dinh2}
\textit{T.C. Dinh}, Enveloppe polynomiale d'un compact de longueur 
finie et chaînes
holomorphes \`a bord rectifiable, \`a paraître dans Acta Mathematica,
\textbf{180:1} (1998).
\bibitem{Dinh3}
\textit{T.C. Dinh}, Orthogonal measures on the boundary of a Riemann
surface and polynomial hull of compacts of finite length,
Pr\'epublication de Paris 6, \textbf{115} (1997), \`a para\^{\i}tre
dans Journal of Functional Analysis.
\bibitem{Dinh4}
\textit{T.C. Dinh}, Probl\`eme du bord dans l'espace projectif
complexe, Pr\'epublication de Paris 6, \textbf{128} (1997).
\bibitem{DolbeaultHenkin}
\textit{P. Dolbeault et G. Henkin}, Chaînes holomorphes de bord donn\'e 
dans $\mathbb{CP}^n$, Bull. Soc. Math. de France, \textbf{125} (1997),
383-445.
\bibitem{Federer}
\textit{F. Federer}, Geometric Measure Theory, Grundlenhren der
Math. Wiss, \textbf{285}, Springer, Berlin-Heidelberg-NewYork, (1988).
\bibitem{GlobevnikStout}
\textit{J. Globevnik, E.L. Stout}, Boundary Morera theorems for
holomorphic functions of several complex variables, Duke Math.J,
\textbf{64}(1991), 571-615. 
\bibitem{Harvey}
\textit{R. Harvey}, Holomorphic chains and their
boundaries, Proc. Symp. Pure Math., \textbf{30}, vol. 1 (1977),
309-382.
\bibitem{HarveyLawson}
\textit{R. Harvey and B. Lawson}, On boundaries of complex analytic 
varieties I, Ann. of Math.,
 \textbf{102} (1975), 233--290.
\bibitem{Henkin}
\textit{G. Henkin}, The Abel-Radon transform and several complex
variables, 
Ann. of Math. Stud., \textbf{7} (1995), 223-275.
\bibitem{KytmanovMyslivets}
\textit{A.M. Kytmanov, S.G. Myslivets}, On a certain boundary analogue of the
Morera theorem, Sibirian Math. J., \textbf{36} (1995), n° 6, 1171-1174.
\bibitem{Rudin}
\textit{W. Rudin}, Function Theory in the Unit Ball of $\mathbb{C}^N$,
Springer, New York, 1980.
\bibitem{Stolzenberg}
\textit{G. Stolzenberg}, Uniform approximation on smooth curves, Acta
Math., \textbf{115} (1966), 185--198.
\bibitem{Stout1}
\textit{ E.L. Stout}, The boundary values of holomorphic functions of
  several complex variables, Duke Math.J., \textbf{44}(1977), 105-108. 
\bibitem{Wermer}
\textit{J. Wermer}, The hull of a curve in $\mathbb{C}^n$, 
Ann. of Math., \textbf{68} (1958),
550--561.
\end{thebibliography}
\end{document}